\documentclass[11pt]{article}
\usepackage{amsthm,amsmath,amssymb,anysize}
\usepackage{graphicx}

\def\qed{\hbox to 0pt{}\hfill$\rlap{$\sqcap$}\sqcup$}

\usepackage[numbers,sort&compress]{natbib}
 \numberwithin{equation}{subsection}

\begin{document}

\title{\textbf{On the prime distribution
  }}

\author{Yong Zhao \quad   Jianqin Zhou\footnote{Corresponding author.  Email: zhou63@ahut.edu.cn}\\
\\
\textit{Department of Computer Science, Anhui University of   Technology}\\
\textit{ Ma'anshan 243002, P. R. China}}

\date{ }
\maketitle
\begin{quotation}
\small\noindent

  In this paper, the estimation formula of the number of primes in a given interval is obtained by using the prime distribution property. For any prime pairs $p>5$ and $ q>5 $, construct a disjoint infinite set sequence $A_1, A_2, \ldots, A_i. \ldots $, such that the number of prime pairs ($p_i$ and $q_i $, $p_i-q_i = p-q $) in $A_i $ increases gradually, where $i>0$. So twin prime conjecture is true. We also prove that for any even integer   $m>2700$,  there exist   more than 10
 prime pairs $(p,q)$, such that $p+q=m$. Thus Goldbach conjecture is true.

\



\end{quotation}

\setcounter{section}{0}\section{Introduction}

Like Goldbach conjecture, twin prime conjecture is also one of the famous unsolved problems in number theory. In 1973,
Chen \cite{Chen}   proved that for any even number $h$, there are  infinite  prime numbers $p$, so that the number of prime factors of $p+h$ does not exceed 2.
In 2008,  Green and Tao \cite{Green}    proved the existence of arbitrarily long arithmetic progressions in the primes.
In 2014, Zhang \cite{Zhang}  proved that bounded gaps
between primes are all less than 70 million.

In this paper, the estimation formula of the number of primes in a given interval is given by using the prime distribution property. For any prime pairs $p>5$ and $q>5$, construct a disjoint infinite set sequence $A_1, A_2, \ldots, A_i. \ldots $, such that the number of prime pairs ($p_i$ and $\ q_i $, $p_i-q_i = p-q $) in $A_i $ increases gradually, where $i>0$. So the original conjecture is true. We also prove that for any even integer   $m>2700$,  there exist   more than 10
 prime pairs $(p,q)$, such that $p+q=m$. Thus Goldbach conjecture is true. For any triple primes $p>7,$  $ q>7 $ and $r>7$, construct a disjoint infinite set sequence $A_1, A_2, \ldots, A_i. \ldots $, such that the number of triple primes ($p_i$  $q_i $ and $r_i $, $p_i-q_i = p-q $,  $p_i-r_i = p-r $) in $A_i $ increases gradually, where $i>0$. So triple primes conjecture is true.


\section{Three lemmas}

In this section, three lemmas are proved.

\textbf{ Lemma 1}.
Given two coprime natural numbers $p$ and $ q $. If the remainder of natural numbers in the set $A$ with respect to $p $ is evenly distributed, then
the remainder of natural numbers in the set $\{aq+c|a\in A\}$ is still evenly distributed, where $c\geq 0 $ is an integer.
 \begin{proof}Prove by a  contradiction.

 Without loss of generality, suppose that $a_i\equiv i \mod p, a_i\in A$, $0\leq i< p$,  but  the remainder of natural numbers in the set $\{a_i q+c|a_i\in A, 0\leq i< p \}$ is not evenly distributed. Let's say that $a_iq+c$ and $a_jq+c$ have the same remainder about $p$,

 $\Longrightarrow a_iq+c- (a_jq+c)= kp, k $ is an integer. $\Longrightarrow$ $(a_i-a_j)q= kp$.

 As $p$ and $ q $ are coprime natural numbers,  $\Longrightarrow$   $ a_i\equiv a_j \mod p $. This is a contradiction.
 \end{proof}

 Based on Lemma 1, we prove the following Lemma 2.

 \textbf{ Lemma 2}. Given natural numbers $\alpha_1, \alpha_2, \ldots , \alpha_n$ pairwise prime, and  the remainder of natural numbers in the set $A$ with respect to $\alpha_i, 1\le i\le n,$  are evenly distributed, where $|A|=\alpha_1\times \alpha_2\times \ldots \times \alpha_n$. We conclude that \[|\{a| a\in A\mbox{ and }a \not \equiv 0\mod \alpha_i, 1\le i\le n\}|=(\alpha_1-1)\times (\alpha_2-1)\times \ldots \times (\alpha_n-1)\]
  \begin{proof}
  If $n=1$, then $|A|=\alpha_1$.  Obviously, $|\{a| a\in A\mbox{ and }a \not\equiv 0\mod \alpha_1\}|=\alpha_1-1$.

  If $n=2$, then $|A|=\alpha_1 \times \alpha_2$.

  Without loss of generality, suppose that

  \[ A=\left(\begin{array}{ccccc}
1 & 2 & \cdots  & \alpha_1-1 &  \alpha_1\\
\alpha_1+1 & \alpha_1+2 & \cdots  & 2\alpha_1-1 &  2\alpha_1\\
\vdots & \vdots & \ddots  & \vdots &  \vdots\\

(\alpha_2 -2)\alpha_1+1 & (\alpha_2 -2)\alpha_1+2 & \cdots  & (\alpha_2 -1)\alpha_1-1 &  (\alpha_2-1)\alpha_1\\
(\alpha_2-1)\alpha_1+1 & (\alpha_2-1)\alpha_1+2 & \cdots  & \alpha_2\alpha_1-1 &  \alpha_2\alpha_1\\
   \end{array}\right)\]
$\Longrightarrow$

{\small
  \[ \{a| a\in A\mbox{ and }a \not\equiv 0\mod \alpha_1\}=\left(\begin{array}{cccc}
1 & 2 & \cdots  & \alpha_1-1  \\
\alpha_1+1 & \alpha_1+2 & \cdots  & 2\alpha_1-1  \\
\vdots & \vdots & \ddots  & \vdots     \\

(\alpha_2 -2)\alpha_1+1 & (\alpha_2 -2)\alpha_1+2 & \cdots  & (\alpha_2 -1)\alpha_1-1  \\
(\alpha_2-1)\alpha_1+1 & (\alpha_2-1)\alpha_1+2 & \cdots  & \alpha_2\alpha_1-1 \\
   \end{array}\right)\]
   }

 As $\alpha_1$ and $ \alpha_2 $ are coprime natural numbers, and by Lemma 1, $\{i\alpha_1+j| 0\le i<\alpha_2\}\equiv\{0,1,2,\cdots,\alpha_2-1 \} \mod\alpha_2 $, where $ 0<j<\alpha_1 $.

 $\Longrightarrow |\{a| a\in A\mbox{ and }a \not \equiv 0\mod \alpha_i, 1\le i\le 2\}|=(\alpha_1-1)\times (\alpha_2-1)$.

 \

  If $n=k+1$, then $|A|=\alpha_1 \times \alpha_2\times \ldots \times \alpha_k \times \alpha_{k+1} $.

  As $\alpha_1, \alpha_2, \ldots , \alpha_{k+1}$ are pairwise primes, and by Lemma 1,  the remainder of natural numbers in the set $\{i\times \alpha_{k+1}|  0<  i\le \alpha_1 \times \alpha_2\times \ldots \times \alpha_k \}$ with respect to $\alpha_i, 1\le i\le k $,  is   evenly distributed.

$\Longrightarrow $  the remainder of natural numbers in the set $\{a| a\in A\mbox{ and }a \not\equiv 0\mod \alpha_{k+1}\} $ with respect to $\alpha_i, 1\le i\le k $,  is   evenly distributed and $|\{a| a\in A\mbox{ and }a \not\equiv 0\mod \alpha_{k+1}\}|= \alpha_1 \times \alpha_2\times \ldots \times \alpha_k \times (\alpha_{k+1}-1) $.

Similarly,  as $\alpha_1, \alpha_2, \ldots , \alpha_{k}$ are pairwise primes, and by Lemma 1,  the remainder of natural numbers in the set $\{i\times \alpha_{k}|  0<  i\le \alpha_1 \times \alpha_2\times \ldots \times \alpha_{k-1} \times (\alpha_{k+1}-1)\}$ with respect to $\alpha_i, 1\le i\le k-1 $,  is   evenly distributed.

$\Longrightarrow $  the remainder of natural numbers in the set $\{a| a\in A\mbox{ and }a \not\equiv 0\mod \alpha_{x}, x=k, k+1\} $ with respect to $\alpha_i, 1\le i\le k-1 $,  is   evenly distributed and $|\{a| a\in A\mbox{ and }a \not\equiv 0\mod \alpha_{x}, x=k, k+1\}|= \alpha_1 \times \alpha_2\times \ldots\times \alpha_{k-1} \times (\alpha_k-1) \times (\alpha_{k+1}-1) $.

$\ldots\ldots$

Finally, we have
\[|\{a| a\in A\mbox{ and }a \not \equiv 0\mod \alpha_i, 1\le i\le k+1\}|=(\alpha_1-1)\times (\alpha_2-1)\times \ldots \times (\alpha_{k+1}-1)\]
   \end{proof}

   The following examples is helpful to understand Lemma 2.

It is known that the natural numbers $2,3,5 $ are mutually prime. Let  $A = \{1,2,3,$ $ \ldots, 29,30 \} $. Obviously
the remainder of  the natural numbers in $A$ about $2,3,5 $ is evenly distributed. According to Lemma 2, after removing all natural numbers with 0 remainder about $2,3,5 $, the number of natural numbers in $A $ becomes $(2-1) \times (3-1) \times (5-1) = 8 $. Namely, $\{ 1,7,11,13,17,19,23,29\}$.

Let $A = \{7*i+1| 1 \le i\le 30 \} $. According to Lemma 1, the remainder of natural numbers about $2,3,5 $ in $A $ is still evenly distributed.
Then, according to Lemma 2, after removing all natural numbers with 0 remainder about $2, 3, 5 $, the number of natural numbers in $A$ becomes $(2-1) \times (3-1) \times (5-1) = 8 $.

\

Based on Lemma 2, the number of primes in 30 consecutive natural numbers $\{6,7, \ldots,$ $ 34,35 \}$ is
 $(2-1)(3-1)(5-1)=30(1-\frac{1}{2})(1-\frac{1}{3})(1-\frac{1}{5})=8$. As $30=2\times3\times5$, the formula is accurate.
 These specific prime numbers are: $7, 11,13, 17,19, 23, 29, 31$.

We are concerned about the following two issues here.

1. Number counting formula of primes $(2-1) (3-1) (5-1)$ is valid for 30 consecutive natural numbers less than 49.
Because 49 is not a multiple of $2,3,5$, but $49=7\times7$. For example, the number of primes   among 30 consecutive natural numbers $\{20, 21, \ldots, 48, 49 \} $ is 7.
These specific prime numbers are:
$23, 29,31, 37,  41, 43, 47$. That is, the actual number of  primes is less than $(2-1)(3-1)(5-1)$.

2. Note that number counting formula of primes \[30(1-\frac{1}{2})(1-\frac{1}{3})(1-\frac{1}{5})(1-\frac{1}{7})\approx 6.86\]
is less than the actual  number 8 of primes in 30 consecutive natural numbers $\{6,7, \ldots,$ $ 34,35 \} $, and also
less than the actual  number 7 of primes in 30 consecutive natural numbers
 $\{20, 21, \ldots, 48, 49\}$.

The number counting formula of primes $30(1-\frac{1}{2})(1-\frac{1}{3})(1-\frac{1}{5})(1-\frac{1}{7})$ is valid for 30 consecutive natural numbers less than $11^2$.
However, it is well known that the prime distribution   is not even,  and the prime numbers in the front part of the effective range are dense and the prime numbers behind are sparse. For example, the actual  number 6 of primes in 30 consecutive natural numbers
 $\{71, 72, \ldots, 99, 100\}$ is less than 6.86. These specific prime numbers are:
$71, 73, 79, 83,  89, 97$.

This paper mainly considers the case that \textbf{the estimation formula for the number of primes   lower than the actual number of primes.}

\

The following proves Lemma 3

\textbf{ Lemma 3}.
\[(1-\frac{d}{30\times x +c})<(1-\frac{d}{30\times (2x+e+1) +c})(1-\frac{d}{30\times (2x+e+2) +c})\]
\[x\ge 0, e\ge 0, 0<c<32, 0<d<30x.\]
  \begin{proof}\ First prove: $\frac{d}{30\times (2x+1) +c}+\frac{d}{30\times (2x+2) +c}<\frac{d}{30\times x +c}$
%
%
%
%
%
%
%
%

\[\Longleftrightarrow (b+c)(2b+c+30+2b+c+60)<(2b+c+30)(2b+c+60),  b=30x\]
\[\Longleftrightarrow (b+c)(4b+2c+90)<(2b+c+30)(2b+c+60)\]
\[\Longleftrightarrow 4b^2+6bc+90b+90c+2c^2 <4b^2+4bc+180b+c^2+90c+1800\]

From both sides of the inequality remove $4b^2+4bc+90b+90c+c^2$

The original inequality $\Longleftrightarrow c^2+2bc<90b+1800 \Longleftrightarrow 0< (90-2c)b+1800-c^2$.

As $c<32$, thus \[\frac{d}{30\times (2x+1) +c}+\frac{d}{30\times (2x+2) +c}<\frac{d}{30\times x +c}\]

\[\Longrightarrow\frac{d}{30\times (2x+e+1) +c}+\frac{d}{30\times (2x+e+2) +c}<\frac{d}{30\times x +c}\]

\[\Longrightarrow (1-\frac{d}{30\times (2x+e+1) +c})(1-\frac{d}{30\times (2x+e+2) +c})\]
\[>1-(\frac{d}{30\times (2x+e+1) +c}+\frac{d}{30\times (2x+e+2) +c})>1-\frac{d}{30\times x +c}\]
  \end{proof}

%
%
%
%
\section{Possible form of prime numbers $\{11+30*x | x\ge 0\}$ and $\{13+30*x | x\ge 0\}$} \label{s1}

Introduce some basic properties of prime numbers. All prime numbers are in odd numbers with single digits of 1, 3, 7 and 9 (except 2 and 5). Now let's take a look at the prime number whose single digit is 1. It is easy to find that there are only two possible forms:

$\{11+30*x | x\ge 0\}$  and  $\{31+30*x | x\ge 0\}$

For prime number whose single digit is 3,  there are only two possible forms:

$\{13+30*x | x\ge 0\}$  and  $\{23 +30*x | x\ge 0\}$

For prime number whose single digit is 7,  there are only two possible forms:

$\{7+30*x | x\ge 0\}$  and  $\{17 +30*x | x\ge 0\}$

For prime number whose single digit is 9,  there are only two possible forms:

$\{19+30*x | x\ge 0\}$  and  $\{29 +30*x | x\ge 0\}$

\

Among possible form of prime numbers $\{11+30*x | x\ge 0\}$,  if any $11+30*x$ is not a prime number, then there are only  four possible decomposition forms:

$H_1=[7+30a][23+30b], a\ge 0, b\ge 0 $;
  $H_2=[13+30a][17+30b], a\ge 0, b\ge 0 $;

    $H_3=[11+30a][31 +30b], a\ge 0, b\ge 0 $;
     $H_4=[19+30a][29+30b], a\ge 0, b\ge 0 $.

As $(19+30a)^2=30c+1$,  thus $(19+30a)^2$ is not a possible decomposition form of $11+30*x$.

 \begin{table*}[h]

 \tiny 
            \begin{center}
                    \caption{ The possible form of prime numbers $\{11+30*x | 0\le x <210\}$}
                    \label{fig:decomposition}
                    \begin{tabular}
                        {|c|c| c|  c|  c|  c|c| c|  c|  c|c|c| c|  c|  c|}
                        \hline
                        11  &  41  &  71  &  101  &  131  &  161  &  191  &  221  &  251  &  281  &  311  &  341  &  371  &  401  &  431   \\
                       \hline
    &     &     &     &     &  7*23 &     &  13*17  &     &     &     &  11*31  &  7*53  &     &         \\
\hline
   461  &  491  &  521  &  551  &  581 & 611  &  641  &  671  &  701  &  731  &  761  &  791  &  821  &  851  &  881   \\
\hline
  &  &  &  19*29  &  7*83  &  13*47  &  &  11*61  &  &  17*43  &  &  7*113  &  &  23*37 &\\

\hline
  911  &  941  &  971  &  1001  &  1031  &  1061  &  1091  &  1121  &  1151  &  1181 & 1211  &  1241  &  1271  &  1301  &  1331\\
 \hline
    &  &  &  7  &  &  &  &  19  &  &  &  7  &  17  &  31  &  &  11\\
    \hline
  1361  &  1391  &  1421  &  1451  &  1481  &  1511  &  1541  &  1571  &  1601  &  1631  &  1661  &  1691  &  1721  &  1751  &  1781  \\

   \hline
    &  13  &  7  &  &  &  &  23  &  &  &  7  &  11  &  19  &  &  17  &  13\\
 \hline

     1811  &  1841  &  1871  &  1901  &  1931  &  1961  &  1991  &  2021  &  2051  &  2081  &  2111  &  2141  &  2171  &  2201  &  2231   \\
      \hline
       &  7  &  &  &  &  37  &  11  &  43  &  7  &  &  &  &  13  &  31  &  23\\
     \hline
        2261  &  2291  &  2321  &  2351  &  2381  &  2411  &  2441  &  2471  &  2501  &  2531  &  2561  &  2591  &  2621  &  2651  &  2681   \\
          \hline
          7  &  29  &  11  &  &  &  &  &  7  &  41  &  &  13  &  &  &  11  &  7\\
     \hline
     2711  &  2741  &  2771  &  2801  &  2831  &  2861  &  2891  &  2921  &  2951  &  2981  &  3011  &  3041  &  3071  &  3101  &  3131   \\
       \hline
       &  &  17  &  &  19  &  &  7  &  23  &  13  &  11  &  &  &  37  &  7  &  31\\

      \hline
      3161  &  3191  &  3221  &  3251  &  3281  &  3311  &  3341  &  3371  &  3401  &  3431  &  3461  &  3491  &  3521  &  3551  &  3581   \\
       \hline
        29  &  &  &  &  17  &  7  &  13  &  &  19  &  47  &  &  &  7  &  53  &\\
      \hline
      3611  &  3641  &  3671  &  3701  &  3731  &  3761  &  3791  &  3821  &  3851  &  3881  &  3911  &  3941  &  3971  &  4001  &  4031   \\
      \hline
        23  &  11  &  &  &  7  &  &  17  &  &  &  &  &  7  &  11  &  &  29\\
      \hline
      4061  &  4091  &  4121  &  4151  &  4181  &  4211  &  4241  &  4271  &  4301  &  4331  &  4361  &  4391  &  4421  &  4451  &  4481   \\
        \hline
        31  &  &  13  &  7  &  37  &  &  &  &  11  &  61  &  7  &  &  &  &\\
        \hline
        4511  &  4541  &  4571  &  4601  &  4631  &  4661  &  4691  &  4721  &  4751  &  4781  &  4811  &  4841  &  4871  &  4901  &  4931  \\
          \hline
          13  &  19  &  7  &  43  &  11  &  59  &  &  &  &  7  &  17  &  47  &  &  13  &\\
      \hline
      4961  &  4991  &  5021  &  5051  &  5081  &  5111  &  5141  &  5171  &  5201  &  5231  &  5261  &  5291  &  5321  &  5351  &  5381   \\
        \hline
       11  &  7  &  &  &  &  19  &  53  &  &  7  &  &  &  11  &  17  &  &\\
        \hline
        5411  &  5441  &  5471  &  5501  &  5531  &  5561  &  5591  &  5621  &  5651  &  5681  &  5711  &  5741  &  5771  &  5801  &  5831   \\
          \hline
          7  &  &  &  &  &  67  &  &  7  &  &  13  &  &  &  29  &  &  7\\
      \hline
      5861  &  5891  &  5921  &  5951  &  5981  &  6011  &  6041  &  6071  &  6101  &  6131  &  6161  &  6191  &  6221  &  6251  &  6281   \\
        \hline
        &  43  &  31  &  11  &  &  &  7  &  13  &  &  &  61  &  41  &  &  7  &  11\\
        \hline

                    \end{tabular}
            \end{center}
        \end{table*}

In Table \ref{fig:decomposition}, the number in the cell indicates that the corresponding $11 + 30 * x $ is not a prime number, and the number in the cell is a factor.

For example, $11 + 30 * 33 = 1001 = 11 \times   7 \times 13$;

$11+30*47=1421=29\times 49=7\times203  $.

Based on Lemma 2, we consider a formula for estimating the number of primes. In Table 1,
consider  first 30 natural numbers $\{11+30*x | 0\le x< 30\}$.
As 30 is not the multiple of $7, 11, 13,17,19, 23,29,31$, thus   consider an estimation formula lower than the actual number of primes:
{\small
\[30(1-\frac{1}{7})(1-\frac{1}{11})(1-\frac{1}{13}) (1-\frac{1}{17})
(1-\frac{1}{19})(1-\frac{1}{23})(1-\frac{1}{29})(1-\frac{1}{31})
 \approx17.20\]
}


The actual number of primes in
$\{11+30*x | 0\le x< 30\}$ is 19.

If $a\in\{11+30*x | 0\le x< 30\}$ and $a$ is a composite number, then $a$ must has one factor in $\{7,11,13,17,19,23\}$, so the above estimation formula is lower than the actual number of primes.

As $(a+30)(b+30)=30(30+a+b)+ab$, the above estimation formula  is valid for  $11+30*x<53\times37=1961=11+30\times65$.
So the above estimation formula is applied to the front part of the effective range, which is the reason that the above estimation formula is lower than the actual number of primes.

\

 \begin{table*}[h]

 \tiny 
            \begin{center}
                    \caption{ The possible form of prime numbers $\{13+30*x | 0\le x <210\}$}
                    \label{fig:decomposition13}
                    \begin{tabular}
                        {|c|c| c|  c|  c|  c|c| c|  c|  c|c|c| c|  c|  c|}
                        \hline
                        13  &  43  &  73  &  103  &  133  &  163  &  193  &  223  &  253  &  283  &  313  &  343  &  373  &  403  &  433   \\
                       \hline
   &  &  &  &  7*19  &  &  &  &  11*23  &  &  &  7*49  &  &  13*31  &\\
\hline
   463  &  493  &  523  &  553  &  583  &  613  &  643  &  673  &  703  &  733  &  763  &  793  &  823  &  853  &  883  \\

\hline
  &  17  &  &  7  &  11  &  &  &  &  19  &  &  7  &  13  &  &  &\\

\hline
  913  &  943  &  973  &  1003  &  1033  &  1063  &  1093  &  1123  &  1153  &  1183  &  1213  &  1243  &  1273  &  1303  &  1333   \\

 \hline
    11  &  23  &  7  &  17  &  &  &  &  &  &  7  &  &  11  &  19  &  &  31\\
    \hline
 1363  &  1393  &  1423  &  1453  &  1483  &  1513  &  1543  &  1573  &  1603  &  1633  &  1663  &  1693  &  1723  &  1753  &  1783  \\

   \hline
   29  &  7  &  &  &  &  17  &  &  11  &  7  &  23  &  &  &  &  &\\
 \hline

     1813  &  1843  &  1873  &  1903  &  1933  &  1963  &  1993  &  2023  &  2053  &  2083  &  2113  &  2143  &  2173  &  2203  &  2233   \\

      \hline
       7  &  19  &  &  11  &  &  13  &  &  7  &  &  &  &  &  41  &  &  7\\
     \hline
       2263  &  2293  &  2323  &  2353  &  2383  &  2413  &  2443  &  2473  &  2503  &  2533  &  2563  &  2593  &  2623  &  2653  &  2683  \\

          \hline
         31  &  &  23  &  13  &  &  19  &  7  &  &  &  17  &  11  &  &  43  &  7  &\\
     \hline
     2713  &  2743  &  2773  &  2803  &  2833  &  2863  &  2893  &  2923  &  2953  &  2983  &  3013  &  3043  &  3073  &  3103  &  3133  \\

       \hline
      &  13  &  47  &  &  &  7  &  11  &  37  &  &  19  &  23  &  17  &  7  &  29  &  13\\

      \hline
     3163  &  3193  &  3223  &  3253  &  3283  &  3313  &  3343  &  3373  &  3403  &  3433  &  3463  &  3493  &  3523  &  3553  &  3583   \\

       \hline
        &  31  &  11  &  &  7  &  &  &  &  41  &  &  &  7  &  13  &  11  &\\
      \hline
      3613  &  3643  &  3673  &  3703  &  3733  &  3763  &  3793  &  3823  &  3853  &  3883  &  3913  &  3943  &  3973  &  4003  &  4033  \\

      \hline
      &  &  &  7  &  &  53  &  &  &  &  11  &  7  &  &  29  &  &  37\\
      \hline

       4063  &  4093  &  4123  &  4153  &  4183  &  4213  &  4243  &  4273  &  4303  &  4333  &  4363  &  4393  &  4423  &  4453  &  4483   \\
               \hline
        17  &  &  7  &  &  47  &  11  &  &  &  13  &  7  &  &  23  &  &  61  &\\

        \hline
        4513  &  4543  &  4573  &  4603  &  4633  &  4663  &  4693  &  4723  &  4753  &  4783  &  4813  &  4843  &  4873  &  4903  &  4933  \\
           \hline
          &  7  &  17  &  &  41  &  &  13  &  &  7  &  &  &  29  &  11  &  &\\
      \hline
     4963  &  4993  &  5023  &  5053  &  5083  &  5113  &  5143  &  5173  &  5203  &  5233  &  5263  &  5293  &  5323  &  5353  &  5383   \\
        \hline
       7  &  &  &  31  &  13  &  &  37  &  7  &  11  &  &  19  &  67  &  &  53  &  7 \\

        \hline
        5413  &  5443  &  5473  &  5503  &  5533  &  5563  &  5593  &  5623  &  5653  &  5683  &  5713  &  5743  &  5773  &  5803  &  5833  \\
            \hline
          &  &  13  &  &  11  &  &  7  &  &  &  &  29  &  &  23  &  7  &  19\\

      \hline
      5863  &  5893  &  5923  &  5953  &  5983  &  6013  &  6043  &  6073  &  6103  &  6133  &  6163  &  6193  &  6223  &  6253  &  6283   \\
         \hline
        11  &  71  &  &  &  31  &  7  &  &  &  17  &  &  &  11  &  7  &  13  &  61\\

        \hline

                    \end{tabular}
            \end{center}
        \end{table*}

Among possible form of prime numbers $\{13+30*x | x\ge 0\}$,  if any $13+30*x$ is not a prime number, then there are only  four possible decomposition forms:

$H_5=[7+30a][19+30b], a\ge 0, b\ge 0 $;
  $H_6=[13+30a][31+30b], a\ge 0, b\ge 0 $;

    $H_7=[11+30a][23 +30b], a\ge 0, b\ge 0 $;
     $H_8=[17+30a][29+30b], a\ge 0, b\ge 0 $.

In Table \ref{fig:decomposition13},
consider  first 30 natural numbers $\{11+30*x | 0\le x< 30\}$.
As 30 is not the multiple of $7, 11, 13,17,19, 23,29,31$, thus   consider an estimation formula lower than the actual number of primes:

{\small
\[30(1-\frac{1}{7})(1-\frac{1}{11})(1-\frac{1}{13}) (1-\frac{1}{17})
(1-\frac{1}{19})(1-\frac{1}{23})(1-\frac{1}{29})(1-\frac{1}{31})
 \approx17.20\]
}

The actual number of primes in
$\{13+30*x | 0\le x< 30\}$ is 20.

If $a\in\{13+30*x | 0\le x< 30\}$ and $a$ is a composite number, then $a$ must has one factor in $\{7,11,13,17,19\}$, so the above estimation formula is lower than the actual number of primes.  As $19+30$ is a composite number, the above estimation formula  is valid for $13+30*x<41\times53=2173=13+30\times72$. So the above estimation formula is applied to the front part of the effective range, which is the reason that the above estimation formula is lower than the actual number of primes.

\

\section{Twin prime conjecture}

We further consider possible form of twin prime numbers $\{(11+30*x,13+30*x) |   x\ge0\}$. From Table \ref{fig:decomposition}  and Table \ref{fig:decomposition13},  (6131, 6133) is  twin prime numbers.

Consider  first 30 pairs of natural numbers $\{(11+30*x,13+30*x) | 0\le x< 30\}$.
As 30 is not the multiple of $7, 11, 13,17,19, 23,29,31$, thus   consider an estimation formula lower than the actual number of twin prime numbers $C(30)$:

{\small
\[30(1-\frac{2}{7})(1-\frac{2}{11})(1-\frac{2}{13}) (1-\frac{2}{17})
(1-\frac{2}{19})(1-\frac{2}{23})(1-\frac{2}{29})(1-\frac{2}{31})
 \approx9.31\]
}

The actual number of twin prime numbers in $\{(11+30*x,13+30*x) | 0\le x< 30\}$ is 13.

The above estimation formula can be understood in this way (take prime number 7 as an example): for every 7 consecutive cells, there must be one cell of $11 + 30 * x_1 $ can be divided by 7, and another cell of $13 + 30 * x_2 $ can be divided by 7, where $x_1\ne x_2$. Therefore, one term in the above formula is $\frac{7-2}{7} $.

When $x = 19 $, $11 + 30 * 19 = 581 $ can be divided by 7, and $13 + 30 * 19 = 583 $ can be divided by 11. Therefore, this cell is counted twice, so the estimation formula $C(30) $ is lower than the actual number of twin primes.

The estimation formula for the   number of twin prime numbers in
 $\{(11+30*x,13+30*x) | 30\le x< 90\}$ is $C(60)$:

 {\tiny

\[60(1-\frac{2}{7})(1-\frac{2}{11})(1-\frac{2}{13}) (1-\frac{2}{17})
(1-\frac{2}{19})(1-\frac{2}{23})(1-\frac{2}{29})(1-\frac{2}{31})\]
\[   (1-\frac{2}{37})(1-\frac{2}{67}) (1-\frac{2}{41})(1-\frac{2}{71})(1-\frac{2}{43})(1-\frac{2}{73})(1-\frac{2}{47})(1-\frac{2}{77})(1-\frac{2}{49})(1-\frac{2}{79})
(1-\frac{2}{53})(1-\frac{2}{83})(1-\frac{2}{59})(1-\frac{2}{89})(1-\frac{2}{61})(1-\frac{2}{91})
     \]
     \[\approx 10.72\]
}

The actual   number of twin prime numbers in
 $\{(11+30*x,13+30*x) | 30\le x< 90\}$   is 15.

In the above formula, $77, 49$ may not appear in the formula because they are multiples of 7. For the completeness of the formula, $77, 49$ are still retained, which only make the valuation smaller.

Because $91^2>13+270*30$,
the effective range of the above estimation formula is: $13+30*x\le13+30\times270$, where $270>90$. So the above estimation formula is applied to the front part of the effective range, which is the main reason that the above estimation formula is lower than the actual number of primes.

The actual   number of twin prime numbers in $\{(11+30*x,13+30*x) | 90\le x< 210\}$ is 24.
The estimation formula for the   number of twin prime numbers in $\{(11+30*x,13+30*x) | 90\le x< 210\}$ is $C(120)$:


\[120(1-\frac{2}{7})\cdots (1-\frac{2}{31})\]
\[(1-\frac{2}{37})(1-\frac{2}{67})\cdots (1-\frac{2}{61})(1-\frac{2}{91})\]
\[(1-\frac{2}{97})(1-\frac{2}{127})(1-\frac{2}{157})(1-\frac{2}{187})\cdots  (1-\frac{2}{121})(1-\frac{2}{151})  (1-\frac{2}{181})(1-\frac{2}{211}) \]


The effective range of the above estimation formula is: $211^2>30*70*210>30*450$.

More generally, \[(30*2^k-30)^2-30*(30*2^{k+1}-30)=900*(2^{2k}-2^{k+2}+2)>0            \]
 \[ \Longrightarrow(30*2^k-30)^2>30*(30*2^{k+1}-30)\]
 where $k\ge2$.

 When $k=2$, the inequality means $90^2 >30*210$.

 When $k=3$, the inequality means $210^2 >30*450$.

  When $k=4$, the inequality means $450^2 >30*930$.

  $\cdots\cdots$

Therefore, the above estimation formula is applied to the front part of the effective range, which is the main reason that the above estimation formula is lower than the actual number of primes.

The actual   number of twin prime numbers in $\{(11+30*x,13+30*x) | 210\le x< 450\}$ is 29. The estimation formula for the   number of twin prime numbers in $\{(11+30*x,13+30*x) | 210\le x< 450\}$ is $C(240)$(only the part about  $\{7+30*x|x\ge0\}$ is given here):

\[(1-\frac{2}{7})\]
\[(1-\frac{2}{37})(1-\frac{2}{67})\]
\[(1-\frac{2}{97})(1-\frac{2}{127})(1-\frac{2}{157})(1-\frac{2}{187})\]
\[(1-\frac{2}{217})(1-\frac{2}{247})(1-\frac{2}{277})(1-\frac{2}{307})(1-\frac{2}{337})(1-\frac{2}{367})(1-\frac{2}{397})(1-\frac{2}{427})\]

It is easy to obtain the following (only the part about  $\{7+30*x|x\ge0\}$ is given here):

\[\frac{C(60)}{C(30)}=2(1-\frac{2}{37})(1-\frac{2}{67})\]
\[\frac{C(120)}{C(60)}=2(1-\frac{2}{97})(1-\frac{2}{127})(1-\frac{2}{157})(1-\frac{2}{187})\]
\[\frac{C(240)}{C(120)}=2(1-\frac{2}{217})(1-\frac{2}{247})(1-\frac{2}{277})(1-\frac{2}{307})(1-\frac{2}{337})(1-\frac{2}{367})(1-\frac{2}{397})(1-\frac{2}{427})\]
\[\cdots\cdots\]

By Lemma 3, \[(1-\frac{2}{217})(1-\frac{2}{247})>(1-\frac{2}{97}), (1-\frac{2}{277})(1-\frac{2}{307})>(1-\frac{2}{127}),\cdots\cdots,\] Thus
\[\cdots\cdots>\frac{C(240)}{C(120)}>\frac{C(120)}{C(60)}>\frac{C(60)}{C(30)}\]

The complete
 $\frac{C(60)}{C(30)}$ is as follows,

 {\tiny

\[  2(1-\frac{2}{37})(1-\frac{2}{67})(1-\frac{2}{41})(1-\frac{2}{71})(1-\frac{2}{43})(1-\frac{2}{73})(1-\frac{2}{47})(1-\frac{2}{77})(1-\frac{2}{49})(1-\frac{2}{79})
(1-\frac{2}{53})(1-\frac{2}{83})(1-\frac{2}{59})(1-\frac{2}{89})(1-\frac{2}{61})(1-\frac{2}{91})
     \]

%
     \[ =2\times 0.575> 1\]
}

Therefore
\[\cdots\cdots> C(240)>C(120) >C(60)>C(30)\approx 9.31\]

In fact, the actual number of twin prime numbers in $\{(11+30*x,13+30*x) | 450\le x< 930\}$ is 71, and the actual number of twin prime numbers in
 $\{(11+30*x,13+30*x) | 930\le x< 1890\}$ is $113, \ldots\ldots,$

Similarly, we can always find another larger interval, which has more than 9.31 twin prime numbers. Therefore, there are infinite twin primes.

\

Among possible form of prime numbers $\{17+30*x | x\ge 0\}$, if any $17+30*x$ is
not a prime number, then there are only four possible decomposition forms:

$H_9=[7+30a][11+30b], a\ge 0, b\ge 0 $;  $H_{10}=[13+30a][29+30b], a\ge 0, b\ge 0 $;

    $H_{11}=[17+30a][31 +30b], a\ge 0, b\ge 0 $;     $H_{12}=[19+30a][23+30b], a\ge 0, b\ge 0 $.

Consider an estimation formula lower
than the actual number of primes in $\{17+30*x | 0\le x< 30\}$:

{\small
\[30(1-\frac{1}{7})(1-\frac{1}{11})(1-\frac{1}{13}) (1-\frac{1}{17})
(1-\frac{1}{19})(1-\frac{1}{23})(1-\frac{1}{29})(1-\frac{1}{31})
 \approx17.20\]
}

 By considering  an estimation formula lower than the actual number of prime pairs    in $\{(11+30*x,17+30*x) | x\ge 0\}$, similarly we can prove
that there are infinite  prime pairs   in $\{(11+30*x,17+30*x) | x\ge 0\}$.

To sum up, we get the following theorem.

\textbf{Theorem 1}. For any two prime numbers $p_0>5$ and $q_0>5$, there are infinite  prime pairs   $p_i$ and $q_i, i\geq 1$, such that $p_i-q_i=p_0-q_0$.

\section{Goldbach conjecture}

Very similar to the case of twin prime conjecture, we further consider possible form of   prime pairs $\{(11+30*x,13+30*(n-x)) |   x\ge0\}$.
From Table \ref{fig:decomposition}  and Table \ref{fig:decomposition13},  (131, 6163) is  a pair of primes, where $n=209, x=4$.

Consider    30 pairs of natural numbers $\{(13+30*x,13+30*(29-x)) | 0\le x< 30\}$.
As 30 is not the multiple of $7, 11, 13,17,19, 23,29,31$, thus   consider an estimation formula lower than the actual number of twin prime numbers $C(30)$:

{\small
\[30(1-\frac{2}{7})(1-\frac{2}{11})(1-\frac{2}{13}) (1-\frac{2}{17})
(1-\frac{2}{19})(1-\frac{2}{23})(1-\frac{2}{29})(1-\frac{2}{31})
 \approx9.31\]
}

The actual number of   prime pairs in $\{(11+30*x,13+30*(29-x)) | 0\le x< 30\}$ is 11.

The above estimation formula can be understood in this way (take prime number 7 as an example): for every 7 consecutive cells, there must be one cell of $11 + 30 * x_1 $ can be divided by 7, and another cell of $13 + 30 * (29-x_2) $ can be divided by 7, where $x_1\ne x_2$. Therefore, one term in the above formula is $\frac{7-2}{7} $.

When $x = 11 $, $11 + 30 * 11 = 341 $ can be divided by 11, and $13 + 30 * (29-11) = 553 $ can be divided by 7. Therefore, this cell is counted twice, so the estimation formula $C(30) $ is lower than the actual number of   prime pairs.

The estimation formula for the   number of prime pairs in
 $\{(11+30*x,13+30*(89-x)) | 30\le x< 90\}$ is $C(60)$:

 {\tiny

\[60(1-\frac{2}{7})(1-\frac{2}{11})(1-\frac{2}{13}) (1-\frac{2}{17})
(1-\frac{2}{19})(1-\frac{2}{23})(1-\frac{2}{29})(1-\frac{2}{31})\]
\[   (1-\frac{2}{37})(1-\frac{2}{67}) (1-\frac{2}{41})(1-\frac{2}{71})(1-\frac{2}{43})(1-\frac{2}{73})(1-\frac{2}{47})(1-\frac{2}{77})(1-\frac{2}{49})(1-\frac{2}{79})
(1-\frac{2}{53})(1-\frac{2}{83})(1-\frac{2}{59})(1-\frac{2}{89})(1-\frac{2}{61})(1-\frac{2}{91})
     \]
     \[\approx 10.72\]
}

The actual   number of prime pairs in
 $\{(11+30*x,13+30*(89-x)) | 30\le x< 90\}$   is 16.

In the above formula, $77, 49$ may not appear in the formula because they are multiples of 7. For the completeness of the formula, $77, 49$ are still retained, which only make the valuation smaller.

The effective range of the above estimation formula is: $91^2>30*270>30*210$.

%

The actual   number of prime pairs in $\{(11+30*x,13+30*(209-x)) | 90\le x< 210\}$ is 29.
The estimation formula for the   number of prime pairs in $\{(11+30*x,13+30*(209-x)) | 90\le x< 210\}$ is $C(120)$:

\[120(1-\frac{2}{7})\cdots (1-\frac{2}{31})\]
\[(1-\frac{2}{37})(1-\frac{2}{67})\cdots (1-\frac{2}{61})(1-\frac{2}{91})\]
\[(1-\frac{2}{97})(1-\frac{2}{127})(1-\frac{2}{157})(1-\frac{2}{187})\cdots  (1-\frac{2}{121})(1-\frac{2}{151})  (1-\frac{2}{181})(1-\frac{2}{211}) \]

%


The effective range of the above estimation formula is: $211^2>30*70*210>30*450$.

More generally, \[(30*2^k-30)^2-30*(30*2^{k+1}-30)=900*(2^{2k}-2^{k+2}+2)>0            \]
 \[ \Longrightarrow(30*2^k-30)^2>30*(30*2^{k+1}-30)\]
 where $k\ge2$.

 When $k=2$, the inequality means $90^2 >30*210$.

 When $k=3$, the inequality means $210^2 >30*450$.

Therefore, the above estimation formula is applied to the front part of the effective range, which is the main reason that the above estimation formula is lower than the actual number of primes.

\

The actual   number of prime pairs in $\{(11+30*x,13+30*(449-x)) | 210\le x< 450\}$ is 44. The estimation formula for the   number of prime pairs in $\{(11+30*x,13+30*(449-x)) | 210\le x< 450\}$ is $C(240)$(only the part about  $\{7+30*x|x\ge0\}$ is given here):

\[(1-\frac{2}{7})\]
\[(1-\frac{2}{37})(1-\frac{2}{67})\]
\[(1-\frac{2}{97})(1-\frac{2}{127})(1-\frac{2}{157})(1-\frac{2}{187})\]
\[(1-\frac{2}{217})(1-\frac{2}{247})(1-\frac{2}{277})(1-\frac{2}{307})(1-\frac{2}{337})(1-\frac{2}{367})(1-\frac{2}{397})(1-\frac{2}{427})\]

It is easy to obtain the following (only the part about  $\{7+30*x|x\ge0\}$ is given here):

\[\frac{C(60)}{C(30)}=2(1-\frac{2}{37})(1-\frac{2}{67})\]
\[\frac{C(120)}{C(60)}=2(1-\frac{2}{97})(1-\frac{2}{127})(1-\frac{2}{157})(1-\frac{2}{187})\]
\[\frac{C(240)}{C(120)}=2(1-\frac{2}{217})(1-\frac{2}{247})(1-\frac{2}{277})(1-\frac{2}{307})(1-\frac{2}{337})(1-\frac{2}{367})(1-\frac{2}{397})(1-\frac{2}{427})\]
\[\cdots\cdots\]

By Lemma 3, \[(1-\frac{2}{217})(1-\frac{2}{247})>(1-\frac{2}{97}), (1-\frac{2}{277})(1-\frac{2}{307})>(1-\frac{2}{127}),\cdots\cdots\] Thus
\[\cdots\cdots>\frac{C(240)}{C(120)}>\frac{C(120)}{C(60)}>\frac{C(60)}{C(30)}\]

The complete
 $\frac{C(60)}{C(30)}$ is as follows,

 {\tiny

\[  2(1-\frac{2}{37})(1-\frac{2}{67})(1-\frac{2}{41})(1-\frac{2}{71})(1-\frac{2}{43})(1-\frac{2}{73})(1-\frac{2}{47})(1-\frac{2}{77})(1-\frac{2}{49})(1-\frac{2}{79})
(1-\frac{2}{53})(1-\frac{2}{83})(1-\frac{2}{59})(1-\frac{2}{89})(1-\frac{2}{61})(1-\frac{2}{91})
     \]

%
     \[ =2\times 0.575> 1\]
}

Therefore
\[\cdots\cdots> C(240)>C(120) >C(60)\approx 10.72\]

In fact, the actual number of prime pairs   in $\{(11+30*x,13+30*(929-x)) | 450\le x< 930\}$ is 73, and the actual number of     prime pairs in
 $\{(11+30*x,13+30*(1889-x)) | 930\le x< 1890\}$ is $136, \ldots\ldots$

Similarly, we can always find another larger interval, which has more than 10.72   prime pairs $(p, q)$, such that $p+q=24+30*(30*2^k-30-1)$, where $k>1$. 

For any $24+30*a, a\ge90$,  if $30*2^k-30-1<a<30*2^{k+1}-30-1$, by Lemma 2,
the estimation formula for the   number of prime pairs in $\{(11+30*x,13+30*(a-x)) | 0\le x\le a\}$ is greater than
 $C(30*2^{k-1})\ge C(60)$.

 Since we can also consider the  prime pairs in $\{(7+30*x,17+30*(a-x)) | 0\le x\le a\}$, thus there exist much more than 10
 prime pairs $(p, q)$, such that $p+q=24+30*a$,


\

For example, for $24+30*99$, since  $30*2^2-30-1<99<30*2^{2+1}-30-1$, so the estimation formula for the   number of prime pairs in $\{(11+30*x,13+30*(99-x)) | 0\le x\le 99\}$ is greater than
 $C(30*2^{2-1})=C(60)$.
In fact, the actual number of prime pairs   in  $\{(11+30*x,13+30*(99-x)) | 0\le x\le 99\}$ is 27, and the actual number of prime pairs   in  $\{(7+30*x,17+30*(99-x)) | 0\le x\le 99\}$ is 32.

\

It is easy to show that for $a>0$,

$\{0+30*a=7+30*(a-1)+23  \}$;

$\{2+30*a=13+30*(a-1)+19  \}$;

$\{4+30*a=11+30*(a-1)+23  \}$;

$\{6+30*a=13+30*(a-1)+23  \}$;

$\{8+30*a=19+30*(a-1)+19 \}$;

$\{10+30*a=17+30*(a-1)+23  \}$;

$\{12+30*a=19+30*(a-1)+23  \}$;

$\{14+30*a=13+30*(a-1)+31  \}$;

$\{16+30*a=17+30*(a-1)+29  \}$;

$\{18+30*a=7+30*a+11  \}$;

$\{20+30*a=7+30*a+13  \}$;

$\{22+30*a=11+30*a+11  \}$;


$\{26+30*a=13+30*a+13  \}$;

$\{28+30*a=11+30*a+17  \}$.

\

Very similar to the case of  $\{24+30*a=11+30*a+ 13  \}$, we can prove that
for any $2l+30*a, 0\le l <15, a\ge90$,   there exist   more than 10
 prime pairs $(p,q)$, such that $p+q=2l+30*a$.

To sum up, we get the following theorem.

\textbf{Theorem 2}. For any even integer   $m>30*90$,  there exist   more than 10
 prime pairs $(p,q)$, such that $p+q=m$.

\section{
Triple prime conjecture}

Based on Lemma 2, we consider a formula for estimating the number of primes. Consider  an estimation formula for the number of primes in first 210 natural numbers $\{11+210*x | 0\le x< 210\}$:

{\footnotesize 
\[210(1-\frac{1}{11})(1-\frac{1}{13}) (1-\frac{1}{17})
(1-\frac{1}{19})(1-\frac{1}{23})(1-\frac{1}{29})(1-\frac{1}{31})(1-\frac{1}{37})(1-\frac{1}{41})(1-\frac{1}{43})
\]
\[(1-\frac{1}{47})(1-\frac{1}{53}) (1-\frac{1}{59})
(1-\frac{1}{61})(1-\frac{1}{67})(1-\frac{1}{71})(1-\frac{1}{73})(1-\frac{1}{79})(1-\frac{1}{83})(1-\frac{1}{89})\]
\[(1-\frac{1}{97})(1-\frac{1}{101}) (1-\frac{1}{103})
(1-\frac{1}{107})(1-\frac{1}{109})(1-\frac{1}{113})(1-\frac{1}{127})(1-\frac{1}{131})(1-\frac{1}{137})(1-\frac{1}{139})(1-\frac{1}{149})\]
\[(1-\frac{1}{151})(1-\frac{1}{157})(1-\frac{1}{163})(1-\frac{1}{167})(1-\frac{1}{173})(1-\frac{1}{179})(1-\frac{1}{181})(1-\frac{1}{191})(1-\frac{1}{193})(1-\frac{1}{197})(1-\frac{1}{199})(1-\frac{1}{211})\]
 \[\approx95.00\]
}


The actual number of primes in
$\{11+210*x | 0\le x< 210\}$ is 98.

%
%

 Similarly, consider  an estimation formula  for the number of primes in first 210 natural numbers $\{17+210*x | 0\le x< 210\}$:

{\footnotesize
\[210(1-\frac{1}{11})(1-\frac{1}{13}) (1-\frac{1}{17})
(1-\frac{1}{19})(1-\frac{1}{23})(1-\frac{1}{29})(1-\frac{1}{31})(1-\frac{1}{37})(1-\frac{1}{41})(1-\frac{1}{43})
\]
\[(1-\frac{1}{47})(1-\frac{1}{53}) (1-\frac{1}{59})
(1-\frac{1}{61})(1-\frac{1}{67})(1-\frac{1}{71})(1-\frac{1}{73})(1-\frac{1}{79})(1-\frac{1}{83})(1-\frac{1}{89})\]
\[(1-\frac{1}{97})(1-\frac{1}{101}) (1-\frac{1}{103})
(1-\frac{1}{107})(1-\frac{1}{109})(1-\frac{1}{113})(1-\frac{1}{127})(1-\frac{1}{131})(1-\frac{1}{137})(1-\frac{1}{139})(1-\frac{1}{149})\]
\[(1-\frac{1}{151})(1-\frac{1}{157})(1-\frac{1}{163})(1-\frac{1}{167})(1-\frac{1}{173})(1-\frac{1}{179})(1-\frac{1}{181})(1-\frac{1}{191})(1-\frac{1}{193})(1-\frac{1}{197})(1-\frac{1}{199})(1-\frac{1}{211})\]
 \[\approx95.00\]
}


The actual number of primes in
$\{17+210*x | 0\le x< 210\}$ is 96.

 Further consider  an estimation formula  for the number of primes in first 210 natural numbers $\{23+210*x | 0\le x< 210\}$:

{\footnotesize
\[210(1-\frac{1}{11})(1-\frac{1}{13}) (1-\frac{1}{17})
(1-\frac{1}{19})(1-\frac{1}{23})(1-\frac{1}{29})(1-\frac{1}{31})(1-\frac{1}{37})(1-\frac{1}{41})(1-\frac{1}{43})
\]
\[(1-\frac{1}{47})(1-\frac{1}{53}) (1-\frac{1}{59})
(1-\frac{1}{61})(1-\frac{1}{67})(1-\frac{1}{71})(1-\frac{1}{73})(1-\frac{1}{79})(1-\frac{1}{83})(1-\frac{1}{89})\]
\[(1-\frac{1}{97})(1-\frac{1}{101}) (1-\frac{1}{103})
(1-\frac{1}{107})(1-\frac{1}{109})(1-\frac{1}{113})(1-\frac{1}{127})(1-\frac{1}{131})(1-\frac{1}{137})(1-\frac{1}{139})(1-\frac{1}{149})\]
\[(1-\frac{1}{151})(1-\frac{1}{157})(1-\frac{1}{163})(1-\frac{1}{167})(1-\frac{1}{173})(1-\frac{1}{179})(1-\frac{1}{181})(1-\frac{1}{191})(1-\frac{1}{193})(1-\frac{1}{197})(1-\frac{1}{199})(1-\frac{1}{211})\]
 \[\approx95.00\]
}


The actual number of primes in
$\{23+210*x | 0\le x< 210\}$ is 94.

\

\

We further consider possible form of triple primes   $\{(11+210*x,17+210*x, 23+210*x) |   x\ge0\}$. For example, when $x=100$, 21011, 21017 and 21023 are all primes.

Consider  an estimation formula $C(210)$ lower than the actual number of triple primes in  first 210 triplets of natural numbers $\{(11+210*x,17+210*x, 23+210*x)| 0\le x< 210\}$:

{\footnotesize
\[210(1-\frac{3}{11})(1-\frac{3}{13}) (1-\frac{3}{17})
(1-\frac{3}{19})(1-\frac{3}{23})(1-\frac{3}{29})(1-\frac{3}{31})(1-\frac{3}{37})(1-\frac{3}{41})(1-\frac{3}{43})
\]
\[(1-\frac{3}{47})(1-\frac{3}{53}) (1-\frac{3}{59})
(1-\frac{3}{61})(1-\frac{3}{67})(1-\frac{3}{71})(1-\frac{3}{73})(1-\frac{3}{79})(1-\frac{3}{83})(1-\frac{3}{89})\]
\[(1-\frac{3}{97})(1-\frac{3}{101}) (1-\frac{3}{103})
(1-\frac{3}{107})(1-\frac{3}{109})(1-\frac{3}{113})(1-\frac{3}{127})(1-\frac{3}{131})(1-\frac{3}{137})(1-\frac{3}{139})(1-\frac{3}{149})\]
\[(1-\frac{3}{151})(1-\frac{3}{157})(1-\frac{3}{163})(1-\frac{3}{167})(1-\frac{3}{173})(1-\frac{3}{179})(1-\frac{3}{181})(1-\frac{3}{191})(1-\frac{3}{193})(1-\frac{3}{197})(1-\frac{3}{199})(1-\frac{3}{211})\]
 \[\approx17.46\]
}

The actual number of triple primes   in  $\{(11+210*x,17+210*x, 23+210*x)| 0\le x< 210\}$ is 19.

Based on Lemma 1,
the above estimation formula can be understood in this way (take prime number 11 as an example): for every 11 consecutive cells, there must be one cell of $11 + 210 * x_1 $ can be divided by 11, another cell of $17 + 210 * x_2 $ can be divided by 11 and the third cell of $23 + 210 * x_3 $ can be divided by 11, where $x_1, x_2, x_3$ are distinct numbers. Therefore, one term in the above formula is $\frac{11-3}{11} $.

However, when $x = 23 $, $11 + 210 * 23 = 4841 $ can be divided by 47,  $17 + 210 * 23 = 4847 $ can be divided by 37,  and $23 + 210 * 23 = 4853 $ can be divided by 23. Therefore, this cell is counted three times, so the estimation formula $C(210) $ is lower than the actual number of triple primes.

The estimation formula for the   number of triple prime numbers in
 $\{(11+210*x,17+210*x, 23+210*x)| 210\le x< 630\}$ is $C(420)$: 

 {\small 

\[420 (1-\frac{3}{11})(1-\frac{3}{13}) \cdots (1-\frac{3}{211})\]
\[   (1-\frac{3}{221})(1-\frac{3}{431}) (1-\frac{3}{223})(1-\frac{3}{433})\cdots (1-\frac{3}{421})(1-\frac{3}{631})\]
     \[\approx 17.66\]
}

The actual   number of triple prime numbers in
  $\{(11+210*x,17+210*x, 23+210*x)| 210\le x< 630\}$   is 22.

In the above formula, $221, 247, 253, \cdots \cdots$ may not appear in the formula because $221=13\times 17, 247=13\times 19, 253=23\times11\cdots\cdots $. For the completeness of the formula, $221, 247, 253, \cdots \cdots$ are still retained, which only make the valuation further smaller.

\

The actual   number of triple prime numbers in
  $\{(11+210*x,17+210*x, 23+210*x)| 630\le x< 1470\}$   is 35.

The estimation formula for the   number of triple prime numbers in  $\{(11+210*x,17+210*x, 23+210*x)| 630\le x< 1470\}$ is $C(840)$:


 {\small 

\[840 (1-\frac{3}{11})(1-\frac{3}{13}) \cdots (1-\frac{3}{211})\]
\[   (1-\frac{3}{221})(1-\frac{3}{431}) (1-\frac{3}{223})(1-\frac{3}{433})\cdots (1-\frac{3}{421})(1-\frac{3}{631})\]
\[   (1-\frac{3}{641})(1-\frac{3}{851})(1-\frac{3}{1061})(1-\frac{3}{1271}) \cdots (1-\frac{3}{841})(1-\frac{3}{1051})(1-\frac{3}{1261})(1-\frac{3}{1471})\]
     \[\approx 20.96\]
}

The effective range of the above estimation formula is: $1471^2>210*7*1470>210*3150$.

More generally, \[(210*2^k-210)^2-210*(210*2^{k+1}-210)=210^2*(2^{2k}-2^{k+2}+2)>0            \]
 \[ \Longrightarrow(210*2^k-210)^2>210*(210*2^{k+1}-210)\]
 where $k\ge2$.

 When $k=2$, the inequality means $630^2 >210*1470$.

 When $k=3$, the inequality means $1470^2 >210*3150$.

 When $k=4$, the inequality means $3150^2 >210*6510$.

  $\cdots\cdots$

Therefore, the above estimation formula is applied to the front part of the effective range, which is the main reason that the above estimation formula is lower than the actual number of primes.

\

It is easy to obtain the following:

{\small
\[\frac{C(420)}{C(210)}=2 (1-\frac{3}{221})(1-\frac{3}{431}) (1-\frac{3}{223})(1-\frac{3}{433})\cdots (1-\frac{3}{421})(1-\frac{3}{631})\]
\[\frac{C(840)}{C(420)}=2(1-\frac{3}{641})(1-\frac{3}{851})(1-\frac{3}{1061})(1-\frac{3}{1271}) \cdots (1-\frac{3}{841})(1-\frac{3}{1051})(1-\frac{3}{1261})(1-\frac{3}{1471})\]
\[\frac{C(1680)}{C(840)}=2(1-\frac{3}{1481}) (1-\frac{3}{1691})   (1-\frac{3}{1901})(1-\frac{3}{2111}) \cdots(1-\frac{3}{2521}) (1-\frac{3}{2731})   (1-\frac{3}{2941})(1-\frac{3}{3151})\]
\[\cdots\cdots\]
}

Similar to the proof of Lemma 3, it is easy to prove that, \[(1-\frac{d}{210\times x +c})<(1-\frac{d}{210\times (2x+e+1) +c})(1-\frac{d}{210\times (2x+e+2) +c})\]
\[x\ge 0, e\ge 0, 0<c<212, 0<d<210x.\]

Therefore,
\[(1-\frac{3}{1481}) (1-\frac{3}{1691})>(1-\frac{3}{641}),   (1-\frac{3}{1901})(1-\frac{3}{2111})>(1-\frac{3}{851}),\cdots\cdots,\] Thus
\[\cdots\cdots>\frac{C(1680)}{C(840)}>\frac{C(840)}{C(420)}>\frac{C(420)}{C(210)}\]

Since
 { \small\small%
\[\frac{C(420)}{C(210)}=2 (1-\frac{3}{221})(1-\frac{3}{431}) (1-\frac{3}{223})(1-\frac{3}{433})\cdots (1-\frac{3}{421})(1-\frac{3}{631})\]
     \[ =2\times 0.506> 1\]
}

Therefore
\[\cdots\cdots> C(1680)>C(840) >C(420)>C(210)\approx 17.46\]

Similarly, we can always find another larger interval, which has more than 17.46 triple prime numbers. Therefore, there are infinite triple primes.

\

 Very similarly, we can further consider possible form of triple primes   $\{(17+210*x,23+210*x, 29+210*x) |   x\ge0\}$,
  $\{(31+210*x,37+210*x, 43+210*x) |   x\ge0\}$,  $\{(47+210*x,53+210*x, 59+210*x) |   x\ge0\}, \cdots\cdots, $  and
   prove
that there are infinite triple primes.

To sum up, we get the following theorem.

\textbf{Theorem 3}. For any triple primes   $p>7,$  $ q>7 $ and $r>7$, there are infinite triple primes   $(p_i, q_i, r_i )$, such that $p_i-q_i = p-q $,  $p_i-r_i = p-r $, where $i>0$.

%

\

\textbf{Jianqin Zhou} (1963-), Ph.D., Professor. He graduated from the Department of mathematics of East China Normal University in 1983 and Fudan University in 1989 with a master degree in mathematics; In 2017, he obtained Ph.D. from Department of computing in Curtin University, Australia. His main research fields are theoretical computer science, combinatorial mathematics and algorithms.

In 1989, he proved a conjecture in combinatorics proposed by the famous mathematician $Paul$ $ Erd\ddot{o}s $ (Please see the paper ``A proof of Alavi conjecture on integer partition", \textit{Acta Mathematicae   Sinica}, 1995,  \textbf{38(5)}636-641). Since 1989, he has published more than 120 papers in \textit{Acta Mathematicae   Sinica}, \textit{Designs, Codes and Cryptography},  \textit{Acta Mathematicae Applicatae Sinica},  \textit{Combinatorica} and other academic journals at home and abroad.

%
%
%
%
%


\vspace{0.2cm}
\begin{thebibliography}{aa}

\bibitem{Chen}
J. Chen,
On the representation of a large even integer as the sum of a prime and the product of at most two primes (in Chinese),
\textit{	Science in China, Series A}. 1973, \textbf{(02)}111-128.

\bibitem{Green}
 B. Green, T. Tao, The primes contain arbitrarily long arithmetic progressions, \textit{Annals of Mathematics}, \textbf{ 167 (2)
(2008)}481-547.

\bibitem{Zhang}
 Y. Zhang, Bounded gaps between primes,
\textit{Annals of Mathematics}, \textbf{179 (2014)}1121-1174.



\end{thebibliography}
\end{document}